\documentclass[11pt,openbib]{article}
\usepackage{amssymb,amsmath,amsfonts,theorem}
\newcommand{\R}{\mathbb{R}}

\newcommand{\gl}{\mathop{\mathrm{\! \, gl}}\nolimits}
\newcommand{\Gl}{\mathop{\mathrm{\! \, Gl}}\nolimits}
\newcommand{\oo}{\mathop{\mathrm{\! \, o}}\nolimits}
\newcommand{\OO}{\mathop{\mathrm{\! \, O}}\nolimits}
\newcommand{\SO}{\mathop{\mathrm{\! \, SO}}\nolimits}
\newcommand{\Ad }{\mathop{\mathrm{Ad}}\nolimits}

\newcommand{\trace}{\mathop{\mathrm{tr}}\nolimits}

\newcommand{\setrule}{\, \rule[-4pt]{.5pt}{13pt}\, }

\newcommand{\onehalf}{\mbox{$\frac{\scriptstyle 1}{\scriptstyle 2}\,$}}
\newcommand{\spann}{\mathop{\rm span}\nolimits}

\newcommand{\vvee}{\mbox{$\scriptscriptstyle \vee $}}

\allowdisplaybreaks
\begin{document}

\begin{center}
{\Large \bf Classification of the coadjoint orbits \\ 
\rule{0pt}{14pt} of the generalized Galile group} \\
\mbox{} \vspace{.05in} \\
Richard Cushman\footnotemark
\end{center}
\footnotetext{printed: \today}
\bigskip 
\begin{abstract} Our goal is to find a representative of each orbit of the 
coadjoint action of the generalized Galile group on the dual of 
its Lie algebra. Our line of argument follows that of Cushman and 
van der Kallen, but differs in the details. 
\end{abstract} 

\section{Basic set up}

Let $\widetilde{V} \subseteq V \subseteq V^{\vvee}$ be a chain of real vector spaces of 
dimension $n$, $n+2$, and $n+3$, respectively, with inner products 
$\widetilde{\gamma }$, $\gamma $, and ${\gamma }^{\vvee}$, whose Gram 
matrix $\widetilde{K}$, $K$, and $K^{\vvee}$ with respect to the standard 
basis $\widetilde{\mathfrak{e}} = \{ e_2, \ldots , e_{n+1} \}$, 
$\mathfrak{e} = \{ e_1, e_2, \ldots , e_{n+1}, e_{n+2} \}$, and 
${\mathfrak{e}}^{\vvee} = \{ e_0, e_1, \ldots , e_{n+1}, e_{n+2} \} $ are $\widetilde{K}$, $K=${\tiny $\begin{pmatrix} 0 & 0 & 1 \\ 0 & \widetilde{K} & 0 \\ 1 & 0 & 0 
\end{pmatrix}$}, and $K^{\vvee}=${\tiny $\left( \begin{array}{c|c} 0 & 0 \\ 
\hline \rule{0pt}{7pt} 0 & K \end{array} \right) $}, respectively. We assume that $K^2 =\mathrm{id}$. Then 
$\widetilde{\gamma }$ and $\gamma $ are nondegenerate and $\ker K^{\vvee} = \spann \{ e_0 \}$. \medskip 

Let $(W, \Gamma )$ be either $(V,\gamma )$ or $(\widetilde{V}, \widetilde{\gamma })$. Then 
\begin{displaymath}
\OO (W,\Gamma ) = 
\{ Q \in \Gl (W, \R ) \setrule \, \Gamma (Qw,Qw') = \Gamma (w,w') 
\, \, \mbox{for all $w,w' \in W$} \} 
\end{displaymath} 
is the \emph{orthogonal group} on $(W,\Gamma )$ with Lie algebra 
\begin{displaymath}
\oo (W, \Gamma ) = \{ Z \in \gl (W, \R ) \setrule \, 
\Gamma (Zw,w') +\Gamma (w,Zw') =0 \, \, \mbox{for all $w,w' \in W$} \} .
\end{displaymath} 
Let 
\begin{displaymath}
\OO (V^{\vvee}, {\gamma }^{\vvee}) = 
\left\{ 
P^{\vvee} \in \Gl (V^{\vvee}, \R ) \setrule \, {\gamma }^{\vvee}(P^{\vvee}v^{\vvee}, P^{\vvee}w^{\vvee}) = 
{\gamma }^{\vvee}(v^{\vvee},w^{\vvee})  \right\} 
\end{displaymath}
be the \emph{generalized orthogonal group} on $(V^{\vvee}, {\gamma }^{\vvee})$. 
Fix $v^{\vvee}_0 \in \ker {\gamma }^{\vvee}$. 
Then $P^{\vvee} v^{\vvee}_0 = v^{\vvee}_0$ for every $P^{\vvee} \in \OO (V^{\vvee}, {\gamma }^{\vvee})$. 
The Lie algebra $\oo (V^{\vvee}, {\gamma }^{\vvee} )$ of $\OO (V^{\vvee}, {\gamma }^{\vvee})$ is equal to 
\begin{displaymath} 
\left\{  X^{\vvee} \in \gl (V^{\vvee}, \R) \setrule \, 
\gamma ^{\vvee}(X^{\vvee}v^{\vvee}, w^{\vvee}) + \gamma ^{\vvee}(v^{\vvee},X^{\vvee}w^{\vvee}) =0 
\, \, \& \, \, X^{\vvee}(v^{\vvee}_0) =0 \right\} .
\end{displaymath}

Let $U = \widetilde{V}, V, V^{\vvee}$ and let $\mathcal{K} = \widetilde{K}, K, K^{\vvee}$. For 
$u = \widetilde{v}, v, v^{\vvee} \in U$ let 
$u^{\ast } = (\mathcal{K}u)^T = u^T\mathcal{K} \in U^{\ast }$. For all $u$, $v \in U$ let 
\begin{displaymath}
L_{u, v} = u \otimes v^{\, \ast} - v \otimes u^{\ast }.
\end{displaymath}
Then $L_{u, v} = - L_{v, u}$. So $L_{u,u} =0$. \medskip 

\noindent \textbf{Fact A.} For every $v^{\vvee}$, $w^{\vvee} \in V^{\vvee}$ we have 
$L_{v^{\vvee},w^{\vvee}} \in \oo (V^{\vvee}, {\gamma }^{\vvee})$. \medskip 

Similarly, for every $v$, $w \in V$ we have $L_{v,w} \in \oo (V,\gamma )$ and for every 
$\widetilde{v}$, $\widetilde{w} \in \widetilde{V}$ we have $L_{\widetilde{v}, 
\widetilde{w}} \in \oo (\widetilde{V}, \widetilde{\gamma })$. \medskip 

\noindent \textbf{Fact B.} If $P^{\vvee} \in \OO (V^{\vvee}, {\gamma }^{\vvee})$, then 
$P^{\vvee}L_{v^{\vvee},w^{\vvee}}(P^{\vvee})^{-1} = L_{P^{\vvee}v^{\vvee}, P^{\vvee}w^{\vvee}}$. \medskip 

Similarly, if $P \in \OO (V, \gamma )$, then 
$PL_{v,w}P^{-1} = L_{Pv,Pw}$. Also if $\widetilde{P} \in \OO (\widetilde{V}, 
\widetilde{\gamma })$, then $\widetilde{P}L_{\widetilde{v}, 
\widetilde{w}}(\widetilde{P})^{-1} = L_{\widetilde{P}\widetilde{v}, 
\widetilde{P}\widetilde{w}}$. \medskip  

With respect to the standard basis $\mathfrak{e}^{\vvee}$ of $V^{\vvee}$ every element 
$P^{\vvee}$ of the generalized orthogonal group $\OO (V^{\vvee}, K^{\vvee})$ is 
given by a matrix {\tiny $\begin{pmatrix} 1 & p^{\ast } \\ 
0 & P \end{pmatrix}$}. Here $p \in V$ and $P \in \OO (V,K)$, using the standard basis $\mathfrak{e}$ of 
$V$. \medskip

\noindent \textbf{Claim 1.} The mapping 
\begin{equation}
\OO (V,K) \ltimes V^{\ast } \rightarrow \OO (V^{\vvee},K^{\vvee}): 
(P, p^{\ast }) \mapsto \mbox{{\footnotesize $\left( \begin{array}{c|c} 
1 & p^{\ast } \\ \hline 
0 & P \end{array} \right)$}} = P^{\vvee} 
\label{eq-s1onenw}
\end{equation}
is an isomorphism of Lie groups, where multiplication $\cdot $ on 
$\OO (V,K) \ltimes V^{\ast }$ is given by 
$(P,p^{\ast }) \cdot (\overline{P}, {\overline{p}}^{\ast} ) = 
(P\overline{P}, {p}^{\ast }\overline{P} + 
{\overline{p}}^{\ast })$ and multiplication on $\OO (V^{\vvee},K^{\vvee})$ is 
matrix multiplication. \medskip

The isomorphism given by equation (\ref{eq-s1onenw}) induces the map 
\begin{equation}
\oo (V,K) \ltimes V^{\ast } \rightarrow \oo (V^{\vvee},K^{\vvee}): (X, x^{\ast}) 
\mapsto \mbox{{\footnotesize $\left( \begin{array}{c|c} 
0 & x^{\ast } \\ \hline   0 & X \end{array} \right) $}} = X^{\vvee}. 
\label{eq-s1twonw}
\end{equation}
The map (\ref{eq-s1twonw}) is an isomorphism of Lie algebras. From now on we will identify $\OO (V^{\vvee},K^{\vvee})$ with $\OO (V,K) \ltimes V^{\ast }$ and 
$\oo (V^{\vvee},K^{\vvee})$ with $\oo (V,K) \ltimes V^{\ast }$. \medskip 

With respect to the standard basis $\mathfrak{e}^{\vvee}$ of $(V^{\vvee},K^{\vvee})$ every 
element $X^{\vvee} = (X, x^{\ast })$ of the Lie algebra $\oo (V^{\vvee},K^{\vvee})$ is given by a matrix 
{\tiny $\left( \begin{array}{c|c} 
0 & {x}^{\ast }  \\ \hline 0 & X \\ \end{array} \right) $}, 
where $x =x_1e_1+\widetilde{x}+x_{n+2}e_{n+2} \in V$, 
which implies $x^{\ast }=x_{n+2}e^T_1+{\widetilde{x}}^{\ast } + 
x_1e^T_{n+2}$. Using the standard basis $\mathfrak{e}$ of $V$, the matrix of $X$ is 
{\tiny $\left( \begin{array}{c|r|c}  
\rule{0pt}{8pt} a & -{\widetilde{e}}^{\ast } & 0 \\  \hline
\rule{0pt}{8pt}\widetilde{d} & \widetilde{X} & \widetilde{e} \\ \hline 
\rule{0pt}{8pt} 0 & -{\widetilde{d}}^{\ast} & -a \end{array} \right)$}
$\in \oo (V,K)$ with $a \in \R$; $\widetilde{d}, \widetilde{e} \in V$;
and $\widetilde{X} \in \oo (\widetilde{V}, \widetilde{K})$. So 
\begin{displaymath}
\dim \oo (V^{\vvee},K^{\vvee}) = n+2 + \dim \oo (V,K) = 3n+3 + \dim \oo (\widetilde{V}, 
\widetilde{K}).
\end{displaymath} 

We now compute ${\Ad }_{P^{\vvee}}X^{\vvee}$, where $X^{\vvee} = (X,x^{\ast }) \in 
\oo (V^{\vvee}, K^{\vvee})$ and $P^{\vvee}=(P,p^{\ast }) \in \OO (V^{\vvee}, K^{\vvee})$. By definition 
\begin{align}
{\Ad }_{P^{\vvee}}X^{\vvee} & = P^{\vvee}X^{\vvee}(P^{\vvee})^{-1} = 
\mbox{{\footnotesize $\begin{pmatrix} 0 & p^{\ast}XP^{-1} + 
x^{\ast }P^{-1} \\ 0 & PXP^{-1} \end{pmatrix} $}.} 
\label{eq-oneA}
\end{align}
Since 
\begin{align}
p^{\ast }XP^{-1} & = p^T(KXP^{-1}) = -p^T(X^TKP^{-1}) \notag \\
& = -p^T(X^TP^TK) = -(PXp)^TK = -(PXp)^{\ast } \notag 
\end{align}
and 
\begin{align}
x^{\ast }P^{-1} = x^TKP^{-1} = x^TP^TK = (Px)^TK = (Px)^{\ast }, 
\notag 
\end{align}
equation (\ref{eq-oneA}) becomes 
\begin{equation}
{\Ad }_{P^{\vvee}}X^{\vvee} = {\Ad }_{(P,p^{\ast })}(X,x^{\ast }) = 
\big( PXP^{-1}, -(PXp)^{\ast } +(Px)^{\ast } \big) .
\label{eq-twoA}
\end{equation}
Replacing $P^{\vvee}$ by $(P^{\vvee})^{-1} = (P^{-1},-(Pp)^{\ast })$ in equation 
(\ref{eq-twoA}) we get  
\begin{equation}
{\Ad}_{(P^{\vvee})^{-1}}X^{\vvee} = \big( P^{-1}XP, (P^{-1}XPp)^{\ast } 
+(P^{-1}x)^{\ast } \big) . 
\label{eq-threeA}
\end{equation}

On $\oo (V^{\vvee},K^{\vvee})$ define an inner product $\langle \, \, | \, \, \rangle $ by 
\begin{displaymath}
\langle (X,x^{\ast }) |(\overline{X}, {\overline{x}}^{\ast }) \rangle = 
\onehalf  \trace (X\overline{X}) + \langle x,\overline{x} \rangle ,
\end{displaymath}
where $(X,x^{\ast }), (\overline{X}, {\overline{x}}^{\ast }) \in 
\oo (V^{\vvee},K^{\vvee})$ and $\langle x, \overline{x} \rangle = {\overline{x}}^TKx$. \medskip

\noindent \textbf{Claim 2.} The inner product $\langle \, \, | \, \, \rangle $ on 
$\oo (V^{\vvee},K^{\vvee})$ is nondegenerate. In other words, the real linear map 
\begin{displaymath}
^{\sharp}: \oo (V^{\vvee},K^{\vvee}) \rightarrow 
{\oo (V^{\vvee},K^{\vvee})}^{\ast }:X^{\vvee} \mapsto (X^{\vvee})^{\sharp}  
\end{displaymath}
is bijective. Here $(X^{\vvee})^{\sharp}({\overline{X}}^{\vvee}) = \langle X^{\vvee} | {\overline{X}}^{\vvee} \rangle $ 
for every ${\overline{X}}^{\vvee} \in \oo (V^{\vvee},K^{\vvee})$. \medskip 

\noindent \textbf{Proof}. Suppose that $0 = \langle (X,x^{\ast }) | 
(\overline{X}, {\overline{x}}^{\ast}) \rangle $ for every 
$(\overline{X}, {\overline{x}}^{\ast}) \in \oo (V^{\vvee},K^{\vvee})$. Let 
${\overline{x}}^{\ast } =0$. Then $0 = \onehalf \trace X\overline{X}$ 
for every $\overline{X} \in \oo (V,K)$. Since $K^2=\mathrm{id}$, 
${\overline{X}}^T \in \oo (V,K)$ if $\overline{X} \in \oo (V,K)$. This follows because $0 = {\overline{X}}^TK + K \overline{X}$ implies 
$0 = \overline{X}K + K{\overline{X}}^T$, after multiplying on the right and 
left by $K$ and using the fact that $K^2 = \mathrm{id}$. Therefore 
$0 = ({\overline{X}}^T)^TK + K{\overline{X}}^T$, that is, ${\overline{X}}^T \in 
\oo (V,K)$. Choose $\overline{X} = X^T$. Then $0 =  
\onehalf \trace X\overline{X} = \onehalf \trace XX^T$, which 
implies $X=0$. Now set $\overline{X} =0$. Then $0 = \langle (0,x^{\ast }) | 
(0, {\overline{x}}^{\ast }) \rangle = \langle x,\overline{x} \rangle = {\overline{x}}^TKx$ 
for every $\overline{x} \in V$. But $K$ is nondegenerate. Therefore $x=0$. 
So $\langle \, \, | \, \, \rangle $ is nondegenerate. 
\hfill $\square $ \medskip 

\noindent \textbf{Fact C.} Let $X \in \oo (V,K)$. For 
$v,w \in V$ we have $\langle Xv,w \rangle = \onehalf \trace  (L_{v,w}X)$. \medskip 

\noindent \textbf{Proof}. For every $y \in V$ we have  
\begin{align}
(L_{v,w}X)y & = L_{v,w}(Xy) = (v\otimes w^{\ast })(Xy) - (w\otimes v^{\ast })(Xy)    \notag \\
& = \langle Xy,w \rangle \, v - \langle Xy,v \rangle \, w \notag \\
& = - \langle y,Xw \rangle \, v + \langle y, Xv \rangle \, w  , \quad 
\mbox{since $X \in \oo (V,K)$} \notag \\
& = -(v\otimes (Xw)^{\ast })y +(w\otimes (Xv)^{\ast })y,  \notag 
\end{align}
that is, $L_{v,w}X = -v\otimes (Xw)^{\ast } + w\otimes (Xv)^{\ast }$. Therefore 
\begin{align} 
\trace (L_{v,w}X) & = -\trace (v\otimes (Xw)^{\ast }) + 
\trace (w \otimes (Xv)^{\ast }) \notag \\
& = -\langle v,Xw \rangle + \langle w, Xv \rangle = \langle Xv,w \rangle + 
\langle w, Xv \rangle = 2 \langle Xv, w \rangle .\tag*{$\square $} 
\end{align}

For each $P^{\vvee} \in \OO (V^{\vvee},K^{\vvee})$ let ${\widetilde{\Ad}}_{P^{\vvee}}$ be 
the mapping of $\oo (V^{\vvee},K^{\vvee})$ into itself such that 
\begin{equation}
({\widetilde{\Ad }}_{P^{\vvee}}(X^{\vvee}))^{\sharp} = 
{\Ad }^T_{(P^{\vvee})^{-1}}\big( (X^{\vvee})^{\sharp } \big) 
\label{eq-oneB}
\end{equation}
for every $X^{\vvee} \in \oo (V^{\vvee},K^{\vvee})$. \medskip

\noindent \textbf{Lemma 3.} For $P^{\vvee} = (P,p^{\ast }) \in \OO (V^{\vvee},K^{\vvee})$ and 
$X^{\vvee} =(X,x^{\, \ast }) \in \oo (V^{\vvee},K^{\vvee})$ we have 
\begin{equation}
{\widetilde{\Ad }}_{(P,p^{\ast })}(X,x^{\ast }) = 
\big( PXP^{-1} +L_{Pp,Px}, (Px)^{\ast } \big) . 
\label{eq-twoB}
\end{equation}

\noindent \textbf{Proof}. Note that the definition of 
the mapping $\widetilde{\Ad}$ is equivalent to 
\begin{align}
\langle {\widetilde{\Ad }}_{(P,p)}(X,x^{\ast }) | 
(Y,y^{\ast }) \rangle & = \langle (X,x^{\ast }) | {\Ad }_{(P^{-1}, 
-(Pp)^{\ast })}(Y,y^{\ast }) \rangle , \notag \\
&\mbox{for every $(Y,y^{\ast }) \in \oo (V^{\vvee},K^{\vvee})_{v^{\vvee}_0}$} 
\notag \\
& \hspace{-1in} = \langle (X,x^{\ast }) | 
(P^{-1}YP, (P^{-1}YPp)^{\ast } + 
(P^{-1}y)^{\ast} ) \rangle , \, \, \mbox{by (\ref{eq-threeA})} \notag \\
&\hspace{-1in} = \onehalf \trace \left[ (XP^{-1}Y)P \right] 
+\langle x,P^{-1}YPp \rangle + \langle x, P^{-1}y \rangle  \notag \\
&\mbox{by definition of $\langle \, \, | \, \, \rangle $} 
\notag \\
& \hspace{-1in} = \onehalf \trace \left[ (PXP^{-1})Y \right]  
+ \langle YPp,Px \rangle + \langle Px,y \rangle \notag \\
&\mbox{since $P \in \OO (V,K)$} \notag \\
&\hspace{-1in} = \onehalf \trace \left[ (PXP^{-1}+L_{Pp,Px})Y \right] 
+\langle Px,y \rangle , \, \, \mbox{see fact C} \notag \\
&\hspace{-1in}  = \langle (PXP^{-1}+L_{Pp,Px}, (Px)^{\ast }) | 
(Y,y^{\ast }) \rangle ,
\notag \\
& \mbox{for every $(Y,y^{\ast}) \in \oo (V^{\vvee},K^{\vvee})$.} \notag
\end{align}
Equation (\ref{eq-twoB}) follows because 
$\langle \, \, | \, \, \rangle $ is nondegenerate.  \hfill $\square $

\section{The generalized Galile group}

First we describe the generalized Galile group in 
terms of matrices. With respect to the standard basis $\mathfrak{e}$ 
of $V$ the \emph{generalized Galile group} ${\mathrm{Gal}}_n$ is 
the set of all $(n+2) \times (n+2)$ real matrices of the form 
{\tiny $\begin{pmatrix} 1 & {\widetilde{a}}^{\ast } & b \\
0 & \widetilde{A} & \widetilde{c} \\ 0 & 0 & 1 \end{pmatrix} $}, where 
$b \in \R $; $\widetilde{a}, \widetilde{c} \in \widetilde{V}$; and 
$\widetilde{A} \in \OO (\widetilde{V},\widetilde{K})$. Here 
${\widetilde{a}}^{\ast } = {\widetilde{a}}^{\, T}\widetilde{K}= 
(\widetilde{K}\, \widetilde{a})^T$.\footnote{In the Galile group one has $\widetilde{A} \in 
\SO (3)$.} \medskip 

Let $y^{\vvee}$ be a nonzero ${\gamma }^{\vvee}$-isotropic vector in $V^{\vvee}$, which does not 
lie in $\ker {\gamma }^{\vvee}$. Let ${\OO (V^{\vvee}, {\gamma }^{\vvee} )}_{y^{\vvee}}$ be the 
subgroup of elements $P^{\vvee}$ of $\OO (V^{\vvee}, {\gamma }^{\vvee})$ such 
that $P^{\vvee}y^{\vvee} = y^{\vvee}$. There is a basis 
${\mathfrak{e}}^{\vvee}$ for $V^{\vvee}$ so that $y^{\vvee} =e_{n+2}$ and the matrix of 
${\gamma }^{\vvee}$ is $K^{\vvee}$. Then 
\begin{align}
{\OO (V^{\vvee},K^{\vvee})}_{e_{n+2}} & = 
\{ P^{\vvee} \in \OO (V^{\vvee},K^{\vvee}) \setrule \, P^{\vvee}e_{n+2} = e_{n+2} \} \notag \\
& \hspace{-.95in} = \{ P^{\vvee} = \mbox{{\footnotesize $\left( 
\begin{array}{c|c} 
1 & p^{\ast } \\ \hline  0 & P \end{array} \right) $}} \in 
\OO (V^{\vvee}, K^{\vvee}) \setrule \, 
p^{\ast }(e_{n+2}) =0, \, \, P \in {\OO (V,K)}_{e_{n+2}} \} . \notag 
\end{align}
Here $P \in {\OO (V,K)}_{e_{n+2}}$ if and only if $P\in \OO (V,K)$ and 
$Pe_{n+2} = e_{n+2}$. To describe an element $(P,p^{\ast })$ of 
${\OO (V^{\vvee},K^{\vvee})}_{e_{n+2}}$ more explicitly we first find the matrix of 
$P\in {\OO (V,K)}_{e_{n+2}}$ with respect to the standard basis 
$\mathfrak{e}$ of $V$. Let $P=${\tiny $ \begin{pmatrix} 
a & {\widetilde{b}}^T & c \\ \widetilde{d} & \widetilde{P} & 
\widetilde{e} \\ f & {\widetilde{g}}^T & h \end{pmatrix} $} be 
a real $(n+2) \times (n+2)$ matrix, where $a$,$c$,$f$,$h \in \R$; 
$\widetilde{b}$,$\widetilde{d}$,$\widetilde{e}$, $\widetilde{g} \in 
\widetilde{V}$; and $\widetilde{P} \in \gl (\widetilde{V}, \R )$. 
Since $Pe_{n+2} = e_{n+2}$, we find that $P=${\tiny $\begin{pmatrix} 
a & {\widetilde{b}}^T & 0 \\ \widetilde{d} & \widetilde{P} & 0 \\
f & {\widetilde{g}}^T & 1 \end{pmatrix} $}. But $P \in \OO (V,K)$, 
that is, $K = P^TKP$. This implies that ${\OO (V,K)}_{e_{n+2}}$ is 
the set of all $(n+2)\times (n+2)$ real matrices of the form 
{\tiny $\begin{pmatrix} 1 & 0 & 0 \\ \widetilde{d} & \widetilde{P} & 0 \\
-\onehalf {\widetilde{d}}^T\widetilde{K}\widetilde{d} & 
- {\widetilde{d}}^{\ast }\widetilde{P} & 1 \end{pmatrix} $}, where 
$\widetilde{d} \in \widetilde{V}$ and $\widetilde{P} \in 
\OO (\widetilde{V}, \widetilde{K})$. If $p = p_1e_1+\widetilde{p} 
+p_{n+2}e_{n+2} \in V$, then $p^{\ast } = p_{n+2} e^T_1 + 
{\widetilde{p}}^{\, \ast } + p_1e^T_{n+2}$. But $0 = p^{\ast }(e_{n+2}) = p_1$. Therefore 
\begin{displaymath}
{\OO (V^{\vvee},K^{\vvee})}_{e_{n+2}} = \Big\{   
\mbox{{\tiny $  \left( \begin{array}{c|ccc} 
1 & p_{n+2} & {\widetilde{p}}^{\ast } & 0 \\ \hline 
0 & 1 & 0 & 0 \\ 
0 & \widetilde{d} & \widetilde{P} & 0 \\
0 &-\onehalf {\widetilde{d}}^T\widetilde{K}\widetilde{d} & 
- {\widetilde{d}}^{\ast }\widetilde{P} & 1 
\end{array} \right) $}\rule[-15pt]{.5pt}{35pt} \, \,   
\raisebox{5pt}{\parbox[t]{1.5in}{$p_{n+2} \in \R$; 
$\widetilde{d}$, $\widetilde{p}\in \widetilde{V}$;\\  
\rule{0pt}{12pt}$\widetilde{P} \in \OO (\widetilde{V}, \widetilde{K}) $}} } \hspace{-15pt} \Big\} .    
\end{displaymath}
Since the map
\begin{equation}
\Xi : {\OO (V^{\vvee},K^{\vvee})}_{e_{n+2}} \rightarrow {\mathrm{Gal}}_n: 
\mbox{{\tiny $ \left( \begin{array}{c|ccc} 
1 & p_{n+2} & {\widetilde{p}}^{\ast } & 0 \\ \hline 
0 & 1 & 0 & 0 \\ 
0 & \widetilde{d} & \widetilde{P} & 0 \\
0 &-\onehalf {\widetilde{d}}^T\widetilde{K}\widetilde{d} & 
- {\widetilde{d}}^{\ast }\widetilde{P} & 1 
\end{array} \right) $}} \mapsto 
\mbox{{\tiny $ \begin{pmatrix} 
1 & {\widetilde{p}}^{\ast } & p_{n+2} \\ 
0 & \widetilde{P} & \widetilde{d} \\
0 & 0 & 1 \end{pmatrix} $}} 
\label{eq-threestar}
\end{equation}
is an isomorphism of groups, we may identify ${\OO (V^{\vvee},K^{\vvee})}_{e_{n+2}}$ with 
the generalized Galile group ${\mathrm{Gal}}_n$. \medskip  

The Lie algebra ${\oo (V^{\vvee},K^{\vvee})}_{e_{n+2}}$ of 
${\OO (V^{\vvee},K^{\vvee})}_{e_{n+2}}$ is 
\begin{displaymath}
\{ X^{\vvee} \in \oo (V^{\vvee},K^{\vvee}) \setrule \, X^{\vvee}(e_{n+2}) 
 =0 \}, 
 \end{displaymath}
 which equals 
\begin{displaymath}
 \{ X^{\vvee}=\mbox{{\footnotesize $\left( \begin{array}{c|c} 
0 & x^{\ast } \\ \hline  0 & X \end{array} \right) $}} \in \oo (V^{\vvee},K^{\vvee}) 
\setrule \, x^{\ast }(e_{n+2}) =0, \, \, X \in {\oo (V,K)}_{e_{n+2}} \} . 
\end{displaymath}
Here $X \in {\oo (V,K)}_{e_{n+2}}$ if and only if $X \in \oo (V,K)$ and 
$X(e_{n+2}) =0$. Since $x^{\ast }(e_{n+2}) =0$ and $X(e_{n+2}) =0$, we 
find that every element $X^{\vvee}$ of ${\oo (V^{\vvee},K^{\vvee})}_{e_{n+2}}$ is given by the matrix 
{\tiny $\begin{pmatrix}
0 & x_{n+2} & {\widetilde{x}}^{\ast } & 0 \\ 
0 & 0 & 0 & 0 \\ 
0 & \widetilde{y} & \widetilde{X} & 0 \\ 
0 & 0 & -{\widetilde{y}}^{\ast }& 0  \end{pmatrix}$}, 
where $x_{n+2} \in \R $; $\widetilde{x}, \widetilde{y} \in \widetilde{V}$; 
and $\widetilde{X} \in \oo (\widetilde{V}, \widetilde{K})$. Thus the map 
\begin{displaymath}
\xi = T_{I}\Xi : {\oo (V^{\vvee}, K^{\vvee})}_{e_{n+2}} \rightarrow {\mathrm{gal}}_n: 
\mbox{{\tiny $\left( \begin{array}{cccc} 
0 & x_{n+2} & {\widetilde{x}}^{\ast } & 0 \\ 
0 & 0 & 0 & 0 \\  
0 & \widetilde{y} & \widetilde{X} & 0 \\ 
0 & 0 & -{\widetilde{y}}^{\ast }& 0  \end{array} \right) $}} \mapsto 
\mbox{\tiny $\left( \begin{array}{ccc} 
0 & {\widetilde{x}}^{\ast } & x_{n+2} \\
0 & \widetilde{X} & \widetilde{y} \\
0 & 0 & 0 \end{array} \right) $} 
\end{displaymath}
is an isomorphism of Lie algebras. Whence 
\begin{displaymath}
\dim {\oo (V^{\vvee},K^{\vvee})}_{e_{n+2}} = 1 + 2n + 
\dim \oo (\widetilde{V}, \widetilde{K}). 
\end{displaymath}

\section{Special cotypes and coadjoint orbits}

Using the standard basis $\mathfrak{e}$ for $V$, the \emph{special tuple} $(V,Y,y;K)$ is formed from 
the pair $(V, Y ;K)$ with $Y \in \oo (V, K)$ and a vector $y \in V$. We say that 
two special tuples $(V, Y,y ;K)$ and 
$(V, Y', y'; K)$ are \emph{equivalent} if and only if there is a $P \in {\OO (V,K)}_{e_{n+2}}$, 
two vectors $v,p \in V$ with $p^{\ast }(e_{n+2}) =0$, and a real number 
$v_0$ such that 
\begin{subequations}
\begin{equation}
Y' +L_{v,e_{n+2}} = P(Y+L_{p,y})P^{-1} 
\label{eq-oneD}
\end{equation}
and
\begin{equation}
y' = Py +v_0e_{n+2}. 
\label{eq-twoD}
\end{equation}
\end{subequations}
Being equivalent is an equivalence relation on the set of special tuples. An equivalence class of special tuples is called a \emph{special cotype}. If $y=0$ and $V = \{ 0 \}$, then $\big( \{ 0 \} , 0, 0; (0) \big)$ is a special tuple, which represents the zero special cotype $\mathbf{0}$. \medskip 

The special tuple $(V,Y,y;K)$ is \emph{decomposable} if there are $Y$-invariant, $K$ orthogonal 
subspaces $V_1$ and $V_2$ of $V$, on which $K$ is nondegenerate, such that $V = V_1 \oplus V_2$ 
with $y \in V_1$ and $V_2 \ne 0$. Let $\nabla $ be the special cotype represented by the tuple $(V,Y,y;K)$. 
The special cotype $\nabla $ is the 
\emph{sum} of the special cotype ${\nabla }'$, represented by the tuple $(V_1, Y|V_1, y; K|V_1)$, and the 
type $\Delta $, represented by the pair $(V_2, Y|V_2; K|V_2)$. We write $\nabla = {\nabla }' + \Delta $. 
If the special cotype $\nabla $ cannot be written as the sum of special cotype ${\nabla }'$ and a 
type $\Delta $, then $\nabla $ is an \emph{indecomposable} special cotype. \medskip

The rest of this section is devoted to proving \medskip

\noindent \textbf{Proposition 4.} The map 
\begin{displaymath}
(V,Y, y ;K) \mapsto  (Y^{\vvee})^{\sharp}|{\oo (V^{\vvee},K^{\vvee})}_{e_{n+2}}, 
\end{displaymath}
where $Y^{\vvee} =(Y,y^{\ast }) \in \oo (V^{\vvee},K^{\vvee})$, induces a bijection 
between the special \linebreak 
cotype $\nabla $, represented by the special tuple $(V,Y,y; K)$, and the 
${\OO (V^{\vvee},K^{\vvee})}_{e_{n+2}}$ coadjoint orbit through 
$(Y^{\vvee})^{\sharp}|{\oo (V^{\vvee},K^{\vvee})}_{e_{n+2}}$.  \medskip 

To prove this proposition we need some preliminary information. \medskip

\noindent \textbf{Fact D.} For $v^{\vvee} = v_0\, e_0 +v_1 \, e_1 +
\widetilde{v} + v_{n+2}\, e_{n+2} \in V^{\vvee}$, 
the matrix of the linear map $L_{v^{\vvee},e_{n+2}}\in \oo (V^{\vvee},K^{\vvee})$ with 
respect to the standard basis ${\mathfrak{e}}^{\vvee}$ of $V^{\vvee}$ is 
{\tiny $\left( \begin{array}{c|ccc}  0 & v_0 & 0 & 0 \\ \hline 0 & v_1 & 0 & 0 \\
0 & \widetilde{v} & 0 & 0 \\ 0 & 0 & -{\widetilde{v}}^{\ast } & - v_1 \end{array} 
\right)  = \left( \begin{array}{c|l} 
0 & v_0 e^T_1 \\ \hline 
\rule{0pt}{12pt}&{\large \mbox{$L_{v,e_{n+2}}$}}  
 \end{array} \right) $} with $v=v_1 \, e_1 + \widetilde{v} + v_{n+2}\, e_{n+2}  \in V$. \medskip

\noindent \textbf{Proof}. By fact A the linear map $L_{v^{\vvee},e_{n+2}}$ of 
$V^{\vvee}$ into itself lies in $\oo (V^{\vvee}, K^{\vvee})$. 
The following computation determines the matrix of $L_{v^{\vvee},e_{n+2}}$ 
with respect to the basis ${\mathfrak{e}}^{\vvee}$. We have  
\begin{align}
L_{v^{\vvee},e_{n+2}}(e_0) & = (v^{\vvee}\otimes e^{\ast }_{n+2})(e_0) 
- (e_{n+2} \otimes (v^{\vvee})^{\ast })(e_0) \notag \\
&= (e^T_{n+2}K^{\vvee}e_0)v^{\vvee} -((v^{\vvee})^TK^{\vvee}e_0)e_{n+2} = 0 \notag 
\end{align}
and
\begin{align} 
L_{v^{\vvee},e_{n+2}}(e_1) & = (e^T_{n+2}K^{\vvee}e_1)v^{\vvee} -((v^{\vvee})^TK^{\vvee}e_1)e_{n+2} 
\notag \\ 
& = e^T_{n+2}(e_{n+2})v^{\vvee} -(v^{\vvee})^T(e_{n+2})e_{n+2} \notag \\
& = v^{\vvee}-v_{n+2}e_{n+2} = v_0e_0 +v_1e_1 + \widetilde{v}. \notag 
\end{align}
For $2 \le i \le n+1$ 
\begin{align} 
L_{v^{\vvee},e_{n+2}}(e_i) & =  (e^T_{n+2}K^{\vvee}e_i)v^{\vvee} -((v^{\vvee})^TK^{\vvee}e_i)e_{n+2} 
= -{\widetilde{v}}^{\, \ast }(e_i)e_{n+2}. \notag 
\end{align}
Finally 
\begin{align} 
L_{v^{\vvee},e_{n+2}}(e_{n+2}) & =  (e^T_{n+2}K^{\vvee}e_{n+2})v^{\vvee} -((v^{\vvee})^TK^{\vvee}e_{n+2})e_{n+2} 
= -v_1e_{n+2}. \tag*{$\square$}
\end{align}

Let ${\oo (V^{\vvee},K^{\vvee})}^0_{e_{n+2}}$ be the set of $(Y^{\vvee})^{\sharp} 
\in {\oo (V^{\vvee},K^{\vvee})}^{\ast }$ with $Y^{\vvee}\in 
\oo (V^{\vvee},K^{\vvee})_{e_{n+2}}$ such that $(Y^{\vvee})^{\sharp }(Z^{\vvee}) =0 $ for every 
$Z^{\vvee} \in {\oo (V^{\vvee},K^{\vvee})}_{e_{n+2}}$. \medskip 

\noindent \textbf{Fact E.} We have 
\begin{displaymath}
{\oo (V^{\vvee},K^{\vvee})}^0_{e_{n+2}} = \{ (L_{v^{\vvee}, e_{n+2}})^{\sharp} \in 
{\oo (V^{\vvee},K^{\vvee})}^{\ast } \setrule \, v^{\vvee} = v_0\, e_0 +v_1\, e_1 
+ \widetilde{v} \in V^{\vvee} \} .
\end{displaymath}

\noindent \textbf{Proof}. Let $Y^{\vvee} = (Y,y^{\ast }) \in \oo (V^{\vvee},K^{\vvee})$. 
Suppose that $0 = (Y^{\vvee})^{\sharp}Z^{\vvee}$ for every $Z^{\vvee} = (Z,z^{\ast }) \in 
{\oo (V^{\vvee},K^{\vvee})}_{e_{n+2}}$. In other words, suppose that $(Y^{\vvee})^{\sharp} 
\in \oo (V^{\vvee},K^{\vvee})^0_{e_{n+2}}$. Then 
$0 = \langle (Y,y^{\ast }) | (Z,z^{\ast }) \rangle = \onehalf \trace YZ + \langle y,z \rangle $ 
for every $Z \in \oo (V,K)$ with $Z(e_{n+2}) =0$ and every $z \in V$ with 
$0  = z^{\ast }(e_{n+2}) = z_1$. So 
for every $\widetilde{Z} \in \oo (\widetilde{V}, \widetilde{K})$, every 
$\widetilde{f}\in \widetilde{V}$, and every $z = \widetilde{z} + 
z_{n+2}e_{n+2} \in V$, we have 
\begin{displaymath}
0 = \onehalf \trace \left[ \mbox{{\tiny $\left( 
\begin{array}{crr} a & -{\widetilde{e}}^{\ast } & 0 \\
\widetilde{d} & \widetilde{Y} & \widetilde{e} \\
0 & -{\widetilde{d}}^{\ast } & - a \end{array} \right) \, 
\left( \begin{array}{crr} 0 & 0 & 0 \\
\widetilde{f} & \widetilde{Z} & 0 \\
0 & -{\widetilde{f}}^{\ast } &  0 \end{array} \right) $}} \right] 
+(0,{\widetilde{z}}^{\, T},z_{n+2})\mbox{{\tiny $\left( 
\begin{array}{ccc} 0 & 0 & 1 \\ 0 & \widetilde{K} & 0 \\ 
1 & 0 & 0 \end{array} \right) $}} \, \mbox{{\footnotesize $ 
\begin{pmatrix} y_1 \\ \widetilde{y} \\ y_{n+2} \end{pmatrix} $,}}  
\end{displaymath}
where $Y=${\tiny $\begin{pmatrix} a & -{\widetilde{e}}^{\ast } & 0 \\
\widetilde{d} & \widetilde{Y} & \widetilde{e} \\ 0 & -{\widetilde{d}}^{\ast } & 
-a \end{pmatrix} $} $\in \oo (V,K)$ and $y =y_1e_1 +\widetilde{y} + 
y_{n+2}e_{n+2} \in V$. In other words, 
\begin{align}
0 & = -\onehalf {\widetilde{e}}^{\, \ast }(\widetilde{f}) + \onehalf \trace (\widetilde{Y}\widetilde{Z}) - 
\onehalf \trace (e \otimes {\widetilde{f}}^{\, \ast }) +y_1z_{n+2}  +
{\widetilde{y}}^{\, T}\widetilde{K}\widetilde{z} \notag \\
& = \onehalf \trace (\widetilde{Y}\widetilde{Z}) 
-{\widetilde{e}}^{\, T}\widetilde{K}\widetilde{f} +y_1z_{n+2} 
+{\widetilde{y}}^{\, T}\widetilde{K}\widetilde{z}, 
\label{eq-twoDdagger}
\end{align}
for every $\widetilde{Z} \in \oo (\widetilde{V}, \widetilde{K})$, 
every $\widetilde{f}, \widetilde{z} \in \widetilde{V}$, and every 
$z_{n+2} \in \R $. Setting $\widetilde{f} = \widetilde{z} =0$ and $z_{n+2} =0$ in 
(\ref{eq-twoDdagger}), we obtain $\widetilde{Y}=0$. Setting $\widetilde{z} =0$ and $z_{n+2} =0$ 
in (\ref{eq-twoDdagger}) with $\widetilde{Y}=0$, we get $\widetilde{f}=0$. Repeating this 
process two more times gives $\widetilde{Y} =0$, $\widetilde{e} =0$, 
$y_1 =0$, and $\widetilde{y}=0$, that is, $Y^{\vvee} = (Y,y^{\ast })$, where 
$Y=${\tiny $ \begin{pmatrix} a & 0 & 0 \\ \widetilde{d} & 0 & 0 \\
0 &-{\widetilde{d}}^{\ast } & -a \end{pmatrix} $}$ \in \oo (V,K)$, $a \in \R $; $\widetilde{d} 
\in \widetilde{V}$; and $y = y_{n+2} e_{n+2}$. Therefore 
\begin{align}
{\oo (V^{\vvee},K^{\vvee})}^0_{e_{n+2}} \subseteq \{ (L_{v^{\vvee}, e_{n+2}})^{\sharp}  \in & \notag \\
&\hspace{-.5in} {\oo (V^{\vvee},K^{\vvee})}^{\ast } \setrule \, v^{\vvee} = y_{n+2} e_0 + a\, e_1 
+ \widetilde{d} \in V^{\vvee} \} .
\label{eq-verynew}
\end{align}
Since
\begin{align}
\dim {\oo (V^{\vvee},K^{\vvee})}^0_{e_{n+2}} &  = \dim \oo (V^{\vvee},K^{\vvee}) 
- \dim {\oo (V^{\vvee},K^{\vvee} )}_{e_{n+2}} 
\notag \\
& = 3+3n + \dim \oo (\widetilde{V}, \widetilde{K}) -1 -2n - \dim \oo (\widetilde{V}, 
\widetilde{K}) \notag \\
& = 2+n \notag 
\end{align}
and the dimension of the subspace of ${\oo (V^{\vvee},K^{\vvee})}^{\ast }$ spanned by 
covectors of the form $(L_{v^{\vvee}, e_{n+2}})^{\sharp}$, where $v^{\vvee} = y_{n+2}e_0 
+ae_1 +\widetilde{d}$, is $2+n$, it follows that equality holds in 
(\ref{eq-verynew}). \hfill $\square $ \medskip 

\noindent \textbf{Proof of proposition 4}. Suppose that the 
special tuple $(V, Y', y'; K)$ is equivalent to the 
special tuple $(V, Y,y; K)$. Then for some $P\in {\OO (V,K)}_{e_{n+2}}$, 
some vectors $-v,p \in V$ with 
$p^{\ast }(e_{n+2}) =0$, and some real number $v_0$ we have 
\begin{displaymath}
Y'+L_{-v,e_{n+2}} = P(Y+L_{p,y})P^{-1} \, \, 
\mathrm{and} \, \, y' = Py +v_0e_{n+2}, 
\end{displaymath}
that is, 
\begin{align}
(Y', (y')^{\ast }) & = 
(PYP^{-1} + L_{Pp,Py}, (Py)^{\ast }) + (L_{v,e_{n+2}}, v_0e^{\ast }_{n+2} 
=v_0e^T_1) 
\label{eq-prop4one}
\end{align}
Now $(Y')^{\vvee} = (Y', (y')^{\ast })$ and 
$Y^{\vvee} = (Y,y^{\ast })$ lie in $\oo (V^{\vvee},K^{\vvee})$. Since $P \in {\OO (V,K)}_{e_{n+2}}$ and 
$p^{\ast }(e_{n+2}) =0$ it follows that $P^{\vvee} = (P,p^{\ast }) \in 
{\OO (V^{\vvee},K^{\vvee})}_{e_{n+2}}$. Let $v^{\vvee}= v_0e_0+v \in V^{\vvee}$. 
Then (\ref{eq-prop4one}) is equivalent to  
\begin{align}
(L_{v^{\vvee},e_{n+2}})^{\sharp } & =  
(Y', (y')^{\ast })^{\sharp} - 
({\widetilde{\Ad }}_{(P,p^{\ast })}(Y,y^{\ast }))^{\sharp } \notag \\
& = ((Y')^{\vvee})^{\sharp } -{\Ad }^T_{(P^{\vvee})^{-1}}(Y^{\vvee})^{\sharp }. 
\end{align}
Consequently, we have
\begin{displaymath}
(Y' )^{\vvee})^{\sharp}|{\oo (V^{\vvee},K^{\vvee})}_{e_{n+2}} = 
{\Ad}^T_{(P^{\vvee})^{-1}}((Y^{\vvee})^{\sharp}|{\oo (V^{\vvee},K^{\vvee})}_{e_{n+2}} ,  
\end{displaymath}
because $(L_{v^{\vvee},e_{n+2}})^{\sharp } \in {\oo (V^{\vvee},K^{\vvee})}^0_{e_{n+2}}$ by fact E. 
In other words, $((Y')^{\vvee})^{\sharp}|$ \linebreak 
${\oo (V^{\vvee},K^{\vvee})}_{e_{n+2}}$ 
lies in the ${\OO (V^{\vvee},K^{\vvee})}_{e_{n+2}}$-coadjoint orbit through 
$(Y^{\vvee})^{\sharp}|$ \linebreak 
${\oo (V^{\vvee},K^{\vvee})}_{e_{n+2}}$. As the special tuple $(V,Y', y'; K)$ 
ranges over the set of special tuples, which are equivalent to $(V,Y,y;K)$, the 
matrix $P \in \OO (V,K)$ and the vector $p\in V$ with $p^{\ast }(e_{n+2}) =0$, 
given by the definition of equivalence, vary so that $P^{\vvee} = (P,p^{\ast })$ 
ranges over all of ${\OO (V^{\vvee},K^{\vvee})}_{e_{n+2}}$. Thus we see that the induced mapping 
is well defined. \medskip  

Suppose that for some $(Y')^{\vvee} \in \oo (V^{\vvee},K^{\vvee})$ and 
some $P^{\vvee} \in {\OO (V^{\vvee},K^{\vvee})}_{e_{n+2}}$ we have 
$((Y')^{\vvee})^{\sharp } = 
{\Ad }^T_{(P^{\vvee})^{-1}}(Y^{\vvee})^{\sharp}$ on ${\oo (V^{\vvee},K^{\vvee})}_{e_{n+2}}$.
In other words, we suppose that $(Y')^{\vvee})^{\sharp }|
{\oo (V^{\vvee},K^{\vvee})}_{e_{n+2}}$ lies in the 
${\OO (V^{\vvee},K^{\vvee})}_{e_{n+2}}$-coadjoint 
orbit through $(Y^{\vvee})^{\sharp}|{\oo (V^{\vvee},K^{\vvee})}_{e_{n+2}}$. 
Then for every $Z^{\vvee} \in {\oo (V^{\vvee},K^{\vvee})}_{e_{n+2}}$ we have 
$\big( (Y')^{\vvee})^{\sharp }-({\widetilde{\Ad}}_{P^{\vvee}}Y^{\vvee})^{\sharp} 
\big) Z^{\vvee} =0$, that is, $((Y')^{\vvee})^{\sharp} - 
({\widetilde{\Ad }}_{P^{\vvee}}Y^{\vvee})^{\sharp} \in {\oo (V^{\vvee},K^{\vvee})}^0_{e_{n+2}}$.  
Consequently, by fact $E$ there is a vector $v^{\vvee} = v_0e_0 -v \in V^{\vvee}$ 
such that $((Y')^{\vvee})^{\sharp } - ({\widetilde{\Ad }}_{P^{\vvee}}Y^{\vvee})^{\sharp} = (L_{v^{\vvee},e_{n+2}})^{\sharp}$. In other words, 
\begin{equation}
(Y')^{\vvee} = {\widetilde{\Ad}}_{P^{\vvee}}Y^{\vvee} +L_{v^{\vvee}, e_{n+2}}. 
\label{eq-oneG}
\end{equation}
Write $(Y')^{\vvee} = (Y', (y')^{\, \ast })$, 
$Y^{\vvee} = (Y,y^{\ast })$ with $Y',Y \in \oo (V,K)$ and 
$y',y \in V$. Also write $P^{\vvee} = (P,p^{\ast })$ with 
$P \in {\OO (V,K)}_{e_{n+2}}$ and $p^{\ast }(e_{n+2}) =0$. Now 
$L_{v^{\vvee},e_{n+2}} = (L_{-v,e_{n+2}}, v_0e^T_1 = v_0 e^{\ast }_{n+2})$. Then 
(\ref{eq-oneG}) may be written as 
\begin{equation}
(Y',(y')^{\, \ast }) = 
(PYP^{-1}+L_{Pp,Py}, (Py)^{\ast }) +(L_{-v,e_{n+2}}, v_0 e^{\ast }_{n+2}), 
\label{eq-twoG}
\end{equation}
that is, 
\begin{subequations}
\begin{equation}
Y' +L_{v, e_{n+2}} = P(Y+L_{p,y})P^{-1} 
\label{eq-threeG}
\end{equation}
and $(y')^{\, \ast } = (Py)^{\ast }+v_0 e^{\ast }_{n+2}$, 
or equivalently, 
\begin{equation}
y' = Py + v_0 e_{n+2}.
\label{eq-fourG}
\end{equation}
\end{subequations}
Now (\ref{eq-threeG}) and (\ref{eq-fourG}) state that the 
special tuples $(V,Y,y;K)$ and $(V, Y', y' ; K)$ are equivalent and thus 
determine the same special cotype. Consequently, the induced mapping is injective. Because every  
element of ${\oo (V^{\vvee},K^{\vvee})}^{\ast }_{e_{n+2}}$ 
can be written as $(Y^{\vvee})^{\sharp }|{\oo (V^{\vvee},K^{\vvee})}_{e_{n+2}}$ for some 
$Y^{\vvee} \in \oo (V^{\vvee},K^{\vvee})$, it follows that induced mapping is surjective. 
So the induced mapping is bijective. This proves proposition 4. \hfill $\square $ \medskip 

\section{Classification of equivalent special tuples}

In view of proposition 4 to classify the coadjoint orbits of ${\OO (V^{\vvee},K^{\vvee})}_{e_{n+2}}$ on 
${\oo (V^{\vvee},K^{\vvee})}^{\ast }_{e_{n+2}}$ we need to classify equivalent special tuples. \medskip 

First we bring the special tuple $(V,Y,y;K)$ into a standard form. Write $y= (y_1, \widetilde{y}, y_{n+2})^T$ and note that the matrix of $P \in {\OO (V,K)}_{e_{n+2}}$ with respect to the standard basis $\mathfrak{e}$ is 
{\tiny $\begin{pmatrix} 
1 & 0 & 0 \\ \widetilde{d} & \widetilde{P} & 0 \\ 
-\onehalf {\widetilde{d}}^{\ast}(\widetilde{d}) & 
-{\widetilde{d}}^{\ast}\widetilde{P} & 1 \end{pmatrix} $}, where 
$\widetilde{d} \in \widetilde{V}$ and $\widetilde{P} \in \OO (\widetilde{V}, 
\widetilde{K})$. Let $(V,Y',y';K)$ be a special tuple which is equivalent to $(V,Y,y;K)$. 
Equation (\ref{eq-twoD}) in the definition of equivalence of 
special tuples reads 
\begin{align}
y' & = Py +v_0e_{n+2} = 
\mbox{{\footnotesize $ \begin{pmatrix} y_1 \\  y_1\widetilde{d} + \widetilde{P}\widetilde{y}\\ 
-\onehalf y_1 {\widetilde{d}}^{\ast }(\widetilde{d}) - 
{\widetilde{d}}^{\ast }(\widetilde{P}\widetilde{y}) 
 + y_{n+2} +v_0 \end{pmatrix} $.}} \label{eq-oneJJ} 
 \end{align}  
Thus we have proved \medskip 

\noindent \textbf{Lemma 5.} Using the standard basis $\mathfrak{e}$ for $V$, the first component $y_1$ of 
the vector $y$ in the special tuple $(V,Y,y;K)$ is invariant under equivalence of special tuples. \medskip

We say that $y_1$ the \emph{parameter} of the special cotype ${\nabla }^{y_1}$, represented by 
the special tuple $(V,Y,y;K)$. \medskip 

There are three cases. \medskip 

\noindent \textsc{case} 1. $y_1\ne 0$. By hypothesis $y_1 \ne 0$, so we may choose 
$\widetilde{d}$ in (\ref{eq-oneJJ}) equal to $-y^{-1}_1 \widetilde{P}\widetilde{y}$. Then the second component of 
$y'$ is $0$. The third component of $y'$ in (\ref{eq-oneJJ}) is 
$\onehalf y^{-1}_1({\widetilde{d}}^{\ast }(\widetilde{d})) +y_{n+2} +v_0$. Now choose $v_0 = 
-\onehalf y^{-1}_1({\widetilde{d}}^{\ast }(\widetilde{d})) -y_{n+2}$. So 
$y' = y_1e_1$. Choosing the vectors $v$ and $p$ in the definition of equivalence equal to $0$, we see that the special tuple $(V,Y,y;K)$ with nonzero parameter $y_1$ is equivalent to the \emph{standardized special tuple} 
$(V,Y', y_1e_1; K)$ with $y_1 \ne 0$ and $Y' = PYP^{-1}$, where $P=${\tiny $\begin{pmatrix} 
1 & 0 & 0 \\ 0 & \widetilde{P} & 0 \\ 0 & 0 & 1 \end{pmatrix}$} with $\widetilde{P} \in \OO (V,K)$. 
Two standardized special tuples $(V,Y,y_1e_1; K)$ and $(V,Y', y_1e_1; K)$ are \emph{equivalent} if there is a $P \in {\OO (V,K)}_{e_{n+2}}$ and a real number $v_0$ such that $Y' = PYP^{-1}$ and 
\begin{align}
& y_1e_1 = y_1P(e_1) +v_0e_{n+2} 
 = y_1\mbox{{\footnotesize $\begin{pmatrix} 1 \\ \widetilde{d} \\
-\onehalf {\widetilde{d}}^{\ast }(\widetilde{d}) +v_0 \end{pmatrix} $}.} 
\label{eq-twoK}
\end{align}
Here $P =${\tiny $\begin{pmatrix} 
1 & 0 & 0 \\ \widetilde{d} & \widetilde{P} & 0 \\ 
-\onehalf {\widetilde{d}}^{\ast}(\widetilde{d}) & 
-{\widetilde{d}}^{\ast}\widetilde{P} & 1 \end{pmatrix} $} with  
$\widetilde{d} \in \widetilde{V}$ and $\widetilde{P} \in \OO (\widetilde{V}, 
\widetilde{K})$. From (\ref{eq-twoK}) we deduce that $y_1\widetilde{d} =0 = 
y_1(-\onehalf {\widetilde{d}}^{\ast }(\widetilde{d})+v_0)$. But $y_1 \ne 0$ 
by hypothesis. So $\widetilde{d} =0$ and consequently $v_0=0$. Therefore 
$P=${\tiny $\begin{pmatrix} 1 & 0 & 0 \\ 
0 & \widetilde{P} & 0 \\ 0 & 0 & 1 \end{pmatrix} $} and $Pe_1 =e_1$. 
Thus the standardized special tuples $(V,Y,y_1e_1;K)$ 
and $(V, Y', y_1e_1; K)$ are equivalent if 
there is a $P \in {\OO (V,K)}_{e_{n+2}}$ with $Pe_1 = e_1$ such that 
$Y'   = PYP^{-1}$. An equivalence class of standardized special tuples is the 
\emph{special cotype} ${\nabla}^{y_1}$. \medskip

Let $(V,Y,y_1e_1;K)$ be a standardized special tuple. Suppose that the matrix of $Y \in \oo (V,K)$ with respect to the standard basis $\mathfrak{e}$ is {\tiny $\begin{pmatrix} \widetilde{a} & 
-{\widetilde{e}}^{\ast } & 0 \\
\widetilde{b} & \widetilde{Y} & \widetilde{e} \\ 0 & -{\widetilde{b}}^{\ast } 
& -\widetilde{a} \end{pmatrix} $}, where $\widetilde{a} \in \R $; $\widetilde{b}, \widetilde{e} \in \widetilde{V}$; and 
$\widetilde{Y} \in \oo (\widetilde{V}, \widetilde{K})$. 
Corresponding to $(V,Y,y_1e_1;K)$ is the \emph{associated pair} 
$(\widetilde{V}, \widetilde{Y}; \widetilde{K})$. Let $(V, Y', (y')_1e_1; K)$ 
be another standardized special tuple. Suppose that the matrix of $Y'$ 
with respect to the standard basis $\mathfrak{e}$ is 
{\tiny $\begin{pmatrix} {\widetilde{a}}' & -({\widetilde{e}}')^{\ast } & 0 \\
{\widetilde{b}}' & {\widetilde{Y}}' & {\widetilde{e}}'\\ 0 & 
-({\widetilde{b}}' )^{\ast } & -{\widetilde{a}}' \end{pmatrix} $}, where ${\widetilde{a}}{\, '} \in \R $, 
${\widetilde{b}}{\, '}, {\widetilde{e}}{\, '} \in \widetilde{V}$, and ${\widetilde{Y}}' \in \oo (\widetilde{V}, 
\widetilde{K})$. Then $(\widetilde{V}, {\widetilde{Y}}', \widetilde{K})$ 
is its associated pair. The associated pairs 
$(\widetilde{V}, \widetilde{Y}; \widetilde{K})$ and 
$(\widetilde{V}, {\widetilde{Y}}', \widetilde{K})$ are \emph{equivalent} if there is $\widetilde{P} \in 
\OO (\widetilde{V}, \widetilde{K})$ such that 
${\widetilde{Y}}' = \widetilde{P}\widetilde{Y}(\widetilde{P})^{-1}$. We call an equivalence class of 
associated pairs an \emph{associated type}, which we denote by 
$\widetilde{\Delta }$. We now prove \medskip

\noindent \textbf{Proposition 6.} The special cotype ${\nabla }^{y_1}$ with nonzero 
parameter $y_1$ uniquely determines, and is uniquely determined by, its 
associated type $\widetilde{\Delta }$ and the nonzero real number $y_1$. \medskip

\noindent \textbf{Proof.} Suppose that the standardized special tuples 
$(V,Y,y_1e_1; K)$ and $(V,Y', y_1e_1, K)$ are equivalent. Then there is $P \in {\OO (V,K)}_{e_{n+2}}$ 
with $Pe_1 = e_1$ and two vectors $v,p \in V$ with $p^{\ast }(e_{n+2}) =p_1 = 0$ 
such that 
\begin{equation}
Y' = PYP^{-1} +L_{Pp,e_1} + L_{-v, e_{n+2}}, 
\label{eq-threeJJ}
\end{equation}
see (\ref{eq-oneD}). Since $Pe_1 = e_1$, equation (\ref{eq-twoD}) with $y=y'=y_1e_1$ implies 
the real number $v_0 =0$. Here $P=$
{\tiny $\begin{pmatrix} 1 & 0 & 0 \\ 0 & \widetilde{P} & 0 \\
0 & 0 & 1 \end{pmatrix} $} with $\widetilde{P} \in \OO (\widetilde{V}, 
\widetilde{K})$. Writing out (\ref{eq-threeJJ}) gives
 \begin{align}
\mbox{\tiny $\begin{pmatrix} {\widetilde{a}}' & (-{\widetilde{e}}')^{\, \ast } & 0 \\ 
{\widetilde{b}}'  & {\widetilde{Y}}' & {\widetilde{e}}' \\ 
0 & -({\widetilde{b}}')^{\, \ast} & -{\widetilde{a}}' \end{pmatrix} $ } &  =   
\mbox{{\tiny $\begin{pmatrix} \widetilde{a} - p_{n+2} +v_1& 
-(\widetilde{P}\widetilde{e})^{\ast } -{\widetilde{p}}^{\ast }& 0 \\ 
\widetilde{P}\widetilde{b}+ \widetilde{v} & \widetilde{P}\widetilde{Y}(\widetilde{P})^{-1} & 
\widetilde{P}\widetilde{e} +\widetilde{p} \\ 0 & -({\widetilde{P}\widetilde{b}})^{\ast } -{\widetilde{v}}^{\ast }& 
-\widetilde{a} + p_{n+2} -v_1
\end{pmatrix} $} }
\label{eq-fourJJ}
\end{align}
Here $p=(0, \widetilde{p}, p_{n+2})^T$ and $v = (v_1, \widetilde{v},0)^T$ are 
vectors in $V$ with respect to the basis $\mathfrak{e}$. Consequently, 
${\widetilde{Y}}' = P\widetilde{Y}P^{-1}$. In other words, 
the associated pairs $(\widetilde{V}, \widetilde{Y}; \widetilde{K})$ 
and $(\widetilde{V},{\widetilde{Y}}'; \widetilde{K})$ 
are equivalent and thus determine a unique associated type $\widetilde{\Delta }$. \medskip  

Now suppose that the pairs $(\widetilde{V}, \widetilde{Y}; \widetilde{K})$ 
and $(\widetilde{V},{\widetilde{Y}}'; \widetilde{K})$ and the nonzero real number $y_1$ associated to the standardized special tuples $(V,Y,y_1e_1; K)$ and $(V,Y', y_1e_1; K)$ are equivalent. Then there is 
$\widetilde{P} \in \OO (\widetilde{V}, \widetilde{K})$ such that 
${\widetilde{Y}}' = \widetilde{P}\widetilde{Y}(\widetilde{P})^{-1}$. We want to determine 
$P\in {\OO (V,K)}_{e_{n+2}}$ with $Pe_1 =e_1$ 
and find two vectors $v,p \in V$ with $p_1 = p^{\ast }(e_{n+2}) =0$ so that 
(\ref{eq-fourJJ}) holds. Let $P=${\tiny $\begin{pmatrix} 1 &0 & 0 \\  0& \widetilde{P} & 0 
\\ 0& 0& 1 \end{pmatrix}$}. Then $P \in {\OO (V,K)}_{e_{n+2}}$ and $Pe_1 = e_1$. 
Now pick $\widetilde{p} =-\widetilde{P}\, \widetilde{e}+ {\widetilde{e}}{\, '}$ and set $p_1 = p_{n+2} =0$. This determines the vector $p$ so that $p^{\ast }(e_{n+2}) =0$. Next choose $v_1 = 
{\widetilde{a}}{\, '} - \widetilde{a}$, $\widetilde{v} = {\widetilde{b}}{\, '} -\widetilde{P}\, \widetilde{b}$, and $v_{n+2} =0$. This determines the vector $v$. For equation (\ref{eq-twoD}) to hold, the real number $v_0 =0$. Using 
${\widetilde{Y} }' = \widetilde{P}\widetilde{Y}(\widetilde{P})^{-1}$ together with the above choices we see that (\ref{eq-fourJJ}) holds. In other words, the standardized special tuples $(V,Y,y_1e_1;K)$ and $(V, Y', y_1e_1; K)$ 
are equivalent and thus determine the special cotype ${\nabla }^{y_1}$. 
\hfill $\square$ \medskip 

If the special tuple $(V,Y,y:K)$ is equivalent to a standardized special tuple $(V,Y,y_1e_1; K)$ with 
$y_1\ne 0$, then the special type ${\nabla }^{y_1}$, represented by $(V,Y,y_1e_1; K)$, is the 
sum of the indecomposable cotype ${\nabla }^{y_1}_2(0)$, $y_1 \ne 0$, represented by the special 
tuple $\big( {\spann}_{\R} \{ e_1, e_{n+2} \} , 0, y_1e_1;${\tiny $\begin{pmatrix} 0 & 1 \\ 1 & 0 \end{pmatrix}$}
$\big)$ with $y_1 \ne 0$, and the associated type $\widetilde{\Delta }$, represented by the pair $(\widetilde{V}, \widetilde{K})$. Note that the dimension of ${\nabla }^{y_1}_2(0), y_1 \ne 0$ is $2$. The nonzero real number $y_1$ is a \emph{modulus} for the cotype ${\nabla }^{y_1}_2(0)$, $y_1 \ne 0$. \medskip 

\noindent This completes case 1. \medskip 

\noindent \textsc{case} 2. $y_1=0$, $y=0$. The special tuple is $(V, Y, 0; K)$ is the sum 
of the zero tuple $\big( \{ 0 \} , 0, 0; (0) \big) $ and the pair $(V, Y; K)$. So the cotype 
$\nabla $, represented by $(V,Y,0; K)$ is the sum of the zero cotype $\mathbf{0}$ and 
the type $\Delta $, represented by the pair $(V,Y; K)$. \medskip 

\noindent This completes case 2. \medskip 

\noindent \textsc{case} 3. $y_1 =0$ and $y \ne 0$. \medskip 

\noindent There are two subcases. \medskip  

\noindent \textsc{subcase} A. We say that the special tuple $(V,Y,y=
${\tiny $\begin{pmatrix} 0 \\ \widetilde{y} \\ y_{n+2} \end{pmatrix}$}; $K)$ with \linebreak
parameter $0$ is an \emph{affine special tuple} if $y$ is a nonzero $K$-isotropic vector. An equivalence class of affine special tuples is an \emph{affine special cotype} $\nabla $. Since $\OO (V ,K)$ acts transitively on the 
collection of all nonzero $K$-isotropic vectors in $V$, there is a $P \in \OO (V, K)$ such that 
$Py = e_{n+2}$. Thus the affine special tuple $(V,Y,y;K)$ is equivalent to the 
\emph{standardized affine special tuple} $(V,Y' = PYP^{-1}, $ $e_{n+2};K)$. To see this we choose $p=0$, $v=0$, and $v_0=0$ in the definition of equivalence of special tuples. \medskip 

Now suppose that two standardized affine special tuples 
$(V,Y,e_{n+2};K)$ and $(V,Y', e_{n+2},K)$ are equivalent. Then there is a $P \in {\OO (V,K)}_{e_{n+2}}$, 
two vectors $-v,p \in V$ with $p_1 = p^{\ast }(e_{n+2}) =0$, and a real number $v_0$ such that 
\begin{subequations}
\begin{equation}
Y' +L_{-v,e_{n+2}} = P(Y+L_{p,e_{n+2}})P^{-1} 
\label{eq-fourJK}
\end{equation}
and 
\begin{equation} 
e_{n+2} = Pe_{n+2} +v_0 e_{n+2} .
\label{eq-fourJL}
\end{equation}
\end{subequations}
From the fact that $Pe_{n+2} = e_{n+2}$, we see equation (\ref{eq-fourJL}) is equivalent to  $v_0 =0$. Therefore we need only consider equation in (\ref{eq-fourJK}). We now write equation (\ref{eq-fourJK}) as  
\begin{align}
Y' & = P(Y+L_{P^{-1}v,e_{n+2}}+L_{p,e_{n+2}})P^{-1}, \quad 
\mbox{since $Pe_{n+2} = e_{n+2}$} \notag \\ 
& = P(Y+L_{w,e_{n+2}})P^{-1}, \notag 
\end{align}
where $w = p+P^{-1}v$. Hence if the standardized affine special tuples 
$(V,Y, e_{n+2};$ $K)$ and $(V,Y', e_{n+2};K)$ are equivalent, then the condition 
\begin{displaymath}
(\ast) \qquad \parbox[t]{3.5in}{there is a $P\in {\OO (V,K)}_{e_{n+2}}$ and 
a vector $w \in V$ such that 
\begin{equation}
Y' = P(Y+L_{w,e_{n+2}})P^{-1}
\label{eq-fourJLstar}
\end{equation} }
\end{displaymath}
holds. Conversely, let $(V,Y, e_{n+2};$ $K)$ and $(V,Y', e_{n+2};K)$ be standardized affine special tuples so that the condition ($\ast$) holds. Set $v_0 =0$. Since $P\in {\OO (V,K)}_{e_{n+2}}$, equation (\ref{eq-fourJL}) holds. Write 
\begin{align}
w & = \mbox{{\tiny $\begin{pmatrix} w_1 \\ \widetilde{w} \\ w_{n+2} 
\end{pmatrix}$}} = \mbox{{\tiny $\begin{pmatrix} w_1 \\ 0 \\ 0\\ 
\end{pmatrix}$}} + \mbox{{\tiny $\begin{pmatrix} 0 \\ \widetilde{w} \\ w_{n+2} 
\end{pmatrix}$}} = P^{-1}v+p, \notag 
\end{align}
where $v =P$\raisebox{2pt}{\tiny $\begin{pmatrix} w_1 \\ 0 \\ 0 \end{pmatrix}$}. Then $p_1 = p^{\ast }(e_{n+2}) =0$. Because $P \in {\OO (V,K)}_{e_{n+2}}$ equation 
(\ref{eq-fourJLstar}) reads
\begin{displaymath}
Y'  = P(Y+L_{P^{-1}v,e_{n+2}}+L_{p,e_{n+2}})P^{-1} = 
L_{v,e_{n+2}} + P(Y+L_{p,e_{n+2}})P^{-1}, 
\end{displaymath}
that is, equation in (\ref{eq-fourJK}) holds. Equation (\ref{eq-fourJLstar}) holds becaue 
$P \in {\OO (V,K)}_{e_{n+2}}$ and $v_0 =0$. Therefore we have proved \medskip 

\noindent \textbf{Lemma 7.} The standardized affine special tuples $(V,Y, e_{n+2};K)$ 
and $(V,Y',$ $e_{n+2};K)$ are equivalent if and only if they satisfy condition ($\ast$) and 
the real number $v_0$ is $0$. \medskip 

Following \cite[p.78]{cushman-vanderkallen} we call an equivalence class of standardized affine special tuples 
an \emph{affine cotype}, which we denote by ${\nabla}_{\mathrm{a}}$. If the special tuple $(V,Y,y;K)$ is an 
affine special tuple, which represents the affine cotype 
${\nabla}_{\mathrm{a}}$, is \emph{decomposable}, then ${\nabla}_{\mathrm{a}} = 
({\nabla }_{\mathrm{a}})' + {\Delta }_3$, where $({\nabla }_{\mathrm{a}})'$ is an affine cotype and 
${\Delta }_3$ is a type. ${\nabla}_{\mathrm{a}}$ is \emph{indecomposable} if no such decomposition 
exists. Every affine cotype ${\nabla}_{\mathrm{a}}$ is the sum of an indecomposable 
affine cotype and a sum of indecomposable types. Up to reordering the summands of this 
decomposition are unique. \medskip

\noindent \textbf{Remark.} Indecomposable affine cotypes are classified in 
\cite[p.83]{cushman-vanderkallen}.\medskip 

\noindent This completes subcase A. \medskip 

\noindent \textsc{subcase} B. Let $(V,Y, y=${\tiny $\begin{pmatrix} 0 \\ \widetilde{y} \\ y_{n+2} \end{pmatrix}$}; $K)$ 
be a special tuple with parameter $0$ which is \emph{not affine}, that is, 
$y$ is not a $K$-isotropic vector. In other words, we have $y^TKy \ne 0$. Such a 
tuple is a \emph{nonaffine} special tuple .\medskip 

\noindent \textbf{Lemma 8.} A nonaffine special tuple is independent of equivalence. \medskip 

\noindent \textbf{Proof.} To see this suppose that the special 
tuple $(V, Y', y';K)$ with parameter $0$ is equivalent 
to the special tuple $(V,Y,y;K)$ with parameter $0$. Then there is 
$P\in {\OO (V,K)}_{e_{n+2}}$, two vectors $v$ and $p$ in $V$ with 
$p_1 = p^{\ast }(e_{n+2}) = 0$, and a real number $v_0$ such that 
$Y' +L_{-v,e_{n+2}} =P(Y+L_{p,e_{n+2}})P^{-1}$ and 
$y' = Py +v_0e_{n+2}$ with $(y')_1 =y_1=0$. So  
\begin{align}
(y')^TKy' & = (Py +v_0e_{n+2})^TK(Py+v_0e_{n+2}) \notag \\
& = (y+v_0e_{n+2})^TP^TKP(y+v_0e_{n+2}), \quad \mbox{since $Pe_{n+2}=e_{n+2}$} \notag \\
& = y^TKy + 2v_0 \, y^TKe_{n+2} +(v_0)^2\, e^T_{n+2}Ke_{n+2}, \quad 
\mbox{since $P\in {\OO (V,K)}_{e_{n+2}}$} \notag \\
& = y^TKy, \quad \parbox[t]{3.5in}{since $Ke_{n+2} = e_1$ implies $y^TKe_{n+2} = 
y^T(e_1) =y_1 =0$ and $e^T_{n+2}Ke_{n+2} =e^T_{n+2}(e_1) =0$} \notag \\
& \ne 0, \quad \mbox{by hypothesis}. \notag 
\end{align}
Therefore the special tuple $(V,Y', y';K)$ with parameter 
$0$ is not affine. \hfill $\square $ \medskip  

Because of lemma 8 we may call an equivalence class of nonaffine 
special tuples a \emph{nonaffine special cotype}, which we denote by ${\nabla}_{\mathrm{na}}$. \medskip 

We now classify indecomposable nonaffine special cotypes, see \cite[p.82]{cushman-vanderkallen}. Suppose that $(V,Y,y;K)$ is a nonaffine special tuple. Since the nonzero vector $y$ is not $K$-isotropic, $y^TKy \ne 0$. Let $f_{n+2} =y$. Therefore on $V_2= 
\spann \{ f_{n+2} \}$ with basis ${\mathfrak{f}}_2 = \{ f_{n+2} \}$ the inner product $K$ is nondegenerate. Let $V_1 = {\spann \{ f_{n+2} \} }^{\perp}$. Then $V_1$ is also $K$-nondegenerate. With respect to the basis  
${\mathfrak{f}}_1 =  \{ f_1, f_2, \ldots , f_{n+1} \}$ of $V_1$ the matrix of 
$K|V_1$ is $F$ and with respect to the basis ${\mathfrak{f}}_2$ of $V_2$ the 
matrix of $K|V_2$ is $\varepsilon {\alpha }^2$, where $\alpha \ne 0$ and 
$\varepsilon = \pm 1$. Therefore with respect to the basis $\mathfrak{f} = 
\{ f_1, \ldots , f_{n+1}, {\alpha }^{-1}f_{n+2} \} $ of $V = V_1 \oplus V_2$ the matrix of $K$ is 
{\tiny $\begin{pmatrix} F & 0 \\ 
0 & \varepsilon  \end {pmatrix} $}. \medskip 

Since $Y \in \oo (V,K)$, the matrix of $Y$ with respect to the basis $\mathfrak{f}$ satisfies 
\begin{align}
0 & = Y^TK +KY = \mbox{{\tiny ${ \begin{pmatrix} Y_1 & \varepsilon \widetilde{b} \\ 
{\widetilde{c}}^T & \widetilde{a} 
\end{pmatrix}}^T \, \begin{pmatrix} F & 0 \\ 0 & \varepsilon  
\end{pmatrix} $}} + \mbox{{\tiny $ \begin{pmatrix} F & 0 \\
0 & \varepsilon  \end{pmatrix} \, \begin{pmatrix} Y_1 & \varepsilon \widetilde{b} \\
{\widetilde{c}}^T & \widetilde{a} \end{pmatrix} $}}, \notag \\
& \hspace{.5in}\mbox{where $Y_1 \in \gl (V_1, \R )$; $\widetilde{b},
\widetilde{c} \in V_1$; and $\widetilde{a} \in \R$} 
\notag \\
& = \mbox{{\tiny $\begin{pmatrix} Y^T_1F & \varepsilon \widetilde{c} \\
\rule{0pt}{10pt} \varepsilon {\widetilde{b}}^TF & \varepsilon  \widetilde{a} \end{pmatrix}$}} + 
\mbox{{\tiny $\begin{pmatrix} FY_1 & \varepsilon F\widetilde{b} \\ 
\rule{0pt}{10pt} {\varepsilon } {\widetilde{c}}^T & 
\varepsilon  \widetilde{a} \end{pmatrix} $}}, \notag 
\end{align}
that is, $Y^T_1F + FY_1 =0$, $\widetilde{c} = -F\widetilde{b} \in V_1$, and $\widetilde{a}=0$. In other words, $Y=${\tiny $\begin{pmatrix} Y_1 & \varepsilon \widetilde{b} \\ 
-(F\widetilde{b})^T & 0 \end{pmatrix} $} with $Y^T_1 F 
+FY_1 =0$ and $\widetilde{b}\in V_1$. \medskip 

Let $w = w' +w_{n+2}f_{n+2} \in V$. We now calculate the matrix of $L_{w, f_{n+2}}$ with respect to the basis $\mathfrak{f}$. For $1 \le i \le n+1$ we have
\begin{align}
L_{w,f_{n+2}}(f_i) & = (w \otimes f^{\ast }_{n+2})(f_i) - (f_{n+2}\otimes w^{\ast })(f_i) \notag \\
& = (f^T_{n+2}Kf_i)w - w^T(Kf_i)f_{n+2} \notag \\ 
& = -(w')^T(Ff_i)f_{n+2} = -(w')^{\ast }(f_i)f_{n+2} \notag 
\end{align}
and
\begin{align}
L_{w,f_{n+2}}(f_{n+2}) & = (w\otimes f^{\ast }_{n+2})(f_{n+2}) - 
(f_{n+2} \otimes w^{\ast })(f_{n+2}) \notag \\
& = (f^T_{n+2}Kf_{n+2})w -(w^TKf_{n+2})f_{n+2} \notag \\
& = \varepsilon w - \varepsilon w_{n+2}f_{n+2} 
= \varepsilon w'. \notag 
\end{align}
Consequently, the matrix of $L_{w, f_{n+2}}$ with respect to the basis 
$\mathfrak{f}$ is {\tiny $ \left( \begin{array}{c|c} 
0 & \varepsilon w' \\ \hline
\rule{0pt}{8pt}(w')^{\, \ast } & 0 \end{array} \right) $}. Since $Y = $\raisebox{2pt}{{\tiny $ \begin{pmatrix} Y_1 & 
\varepsilon \widetilde{b} \\
-(F\widetilde{b})^T & 0 \end{pmatrix}$}}$ = 
$\raisebox{2pt}{{\tiny $\begin{pmatrix} Y_1 & 0 \\ 0 & 0 \end{pmatrix}$}}$+ L_{\widetilde{b}, f_{n+2}}$ with 
$\widetilde{b} \in V_1$, the nonaffine special \linebreak 
tuple $(V,Y,y;K)$ is equivalent to the nonaffine special tuple 
$( V, Y' =${\tiny $\begin{pmatrix} Y_1 & 0 \\ 0 & 0 \end{pmatrix} $}, 
$f_{n+2}; K =${\tiny $\begin{pmatrix} F & 0 \\ 0 & \varepsilon  \end{pmatrix} $}$)$. \medskip 

Now the subspace $V_2 = \spann \{ f_{n+2} \} $ is $K$-nondegenerate. 
Because $Y'f_{n+2}$ $=0$, $V_2$ is $Y'$-invariant. Because $Y' \in \oo (V,K)$ the subspace $V_1 = 
{\spann \{ f_{n+2} \} }^{\perp }$ is $Y'$-invariant and 
$K$-nondegenerate. Thus the nonaffine special cotype ${\nabla}_{\mathrm{na}}$, 
represented by the special tuple $(V,Y,y;K)$, is the sum of the nonaffine cotype 
$({\nabla}_{\mathrm{na}})'$, represented by the nonaffine special tuple 
$\big( V_2, 0,f_{n+2};K|V_2 = (\varepsilon )\big)$, and the type ${\Delta }_1$ represented by 
the pair $(V_1,Y_1=Y|V_1; K|V_1 =F)$. Since $y = \alpha f_{n+2}$, $\alpha $ is a \emph{modulus} for 
$({\nabla}_{\mathrm{na}})'$. Changing the sign of $f_{n+2}$, if necessary, we may assume that 
$\alpha >0$. We use the notation ${\nabla}^{\varepsilon}_{\mathrm{na}}$, $\alpha >0$, for 
$({\nabla }_{\mathrm{na}})'$. Clearly ${\nabla}^{\varepsilon}_{\mathrm{na}}$, $\alpha >0$ is an indecomposable nonaffine special cotype. \medskip 

This completes subcase B and the classification of special cotypes. \hfill $\square $

\section{The main result}

This section is devoted to proving \medskip

\noindent \textbf{Theorem 9.} Let $\nabla $ be a nonzero special cotype. If its 
parameter $y_1$ is nonzero, then $\nabla $ is the sum of the indecomposable type 
${\nabla }^{y_1}_2(0), y_1\ne 0$ and the sum indecomposable associated types. If the parameter of $\nabla $ is $0$ and ${\nabla}$ is a nonaffine special cotype ${\nabla}_{\mathrm{na}}$, then $\nabla $ is the sum of the indecomposable nonaffine special cotype ${\nabla }^{\varepsilon}_{\mathrm{na}}$, $\alpha >0$ and the sum of indecomposable types. If $\nabla $ is an affine special cotype 
${\nabla}_{\mathrm{a}}$, then $\nabla $ is the sum of an indecomposable affine cotype and 
a sum of indecomposable types. Each of the above sums is unique up to 
reordering of the summands. \medskip

\noindent \textbf{Proof.} Suppose that $\nabla$ is a nonzero special cotype represented by the \linebreak 
special tuple $(V,Y,y;K)$. If $y_1$, the first component of $y$ with respect to the standard basis $\mathfrak{e}$, is nonzero, then $\nabla$ is the sum of an indecomposable cotype ${\nabla }^{y_1}_2(0), y_1 \ne 0$, represented by 
the standardized tuple $({\spann }_{\R }\{ e_1, e_{n+2},0, y_1 e_1; $ $K|{\spann }_{\R }\{ e_1, e_{n+2} \} )$, and uniquely determined by the associated type 
$\widetilde{\Delta }$, represented by the associated pair $(\widetilde{V}, \widetilde{Y};\widetilde{K})$. 
The type $\widetilde{\Delta }$ is a unique sum (up to reordering) of indecomposable types, 
see \cite{burgoyne-cushman}. The argument stops here. If $y_1$ is zero and $\nabla $ is a nonaffine special cotype ${\nabla }_{\mathrm{na}}$, then ${\nabla}_{\mathrm{na}} = {\nabla }^{\varepsilon}_{\mathrm{na}}, \alpha >0
+ \Delta $, where the type $\Delta$ is a unique sum (up to reordering) of indecomposable types. The 
argument stops here. If $y_1$ is zero and $\nabla$ is an affine special cotype 
${\nabla }_{\mathrm{a}}$, then ${\nabla }_{\mathrm{a}}$ is uniquely determined by its little cotype 
$({\nabla}_{\mathrm{a}})_{\ell }$, see \cite[p.80]{cushman-vanderkallen}. 
$({\nabla}_{\mathrm{a}})_{\ell } $ is either the zero cotype, a nonaffine cotype, or 
is an affine special cotype. If $({\nabla }_{\mathrm{a}})_{\ell }$ is 
the zero cotype, the argument stops. If the second 
alternative holds, applying lemma 8 one obtains that 
$({\nabla}_{\mathrm{a}})_{\ell }$ is a sum of a unique indecomposable 
nonaffine special cotype and a sum of indecomposable types. The argument then stops. 
Otherwise, repeat the above argument using the affine cotype 
$({\nabla}_{\mathrm{a}})_{\ell }$ instead of the 
cotype ${\nabla }_{\mathrm{a}}$, see \cite[lemma 15,p.83]{cushman-vanderkallen}. This iterative procedure eventually stops because $\dim ({\nabla}_{\mathrm{a}})_{\ell } < 
\dim {\nabla}_{\mathrm{a}}$ and eventually $({\nabla}_{\mathrm{a}})_{\ell }$ is the zero cotype and the 
argument stops. Therefore ${\nabla }_{\mathrm{a}}$ is the unique sum (up to reordering) of an indecomposable affine cotype and a sum of indecomposable types. \hfill $\square $ 

\section{Coadjoint orbits of the general Galile group}

In this section we use the theory of section $4$ to classify the coadjoint orbits of 
the general Galile group ${\mathrm{Gal}}_3 = \OO ({\R }^6, K^{\vvee})_{e_5}$. With 
respect to the standard basis ${\mathfrak{e}}^{\vvee}_5 = \{ e_0, e_1, \ldots , e_5 \} $ of 
${\R }^6$ the matrix $K^{\vvee}$ of the inner product ${\gamma }^{\vvee}$ on $V^{\vvee}$ is 
{\tiny $ \left( \begin{array}{c|c} 
0 &  \\ \hline 
\rule{0pt}{6pt} & K \end{array} \right) =$}{\tiny $\left( \begin{array}{c|ccc} 
  0 & & &  \\ \hline
   &0 &0 & 1 \\  
   & 0 & I_3 & 0 \\
   & 0 & 0 & 1  \end{array} \right) $}, which has signature pair $(p,n) = (4,1)$, nullity $1$, 
index $1$, and dimension $6$. The kernel of ${\gamma }^{\vvee}$ is spanned by $e_0$. \medskip 

In table 1 we list the possible $\oo ({\R }^5, K)$ indecomposable types.
\begin{center}
\begin{tabular}{rccc}
&\multicolumn{1}{c}{type} & \multicolumn{1}{c}{dimension} & \multicolumn{1}{c}{index} \\ \hline
\rule{0pt}{11pt}1. & ${\Delta }^{-}_2(0)$ & $3$ & $1$ \\
2. & ${\Delta}_0(\zeta , \mathrm{RP})$ & $2$ & $1$ \\
3. & ${\Delta }_0(\mathrm{i}\beta , \mathrm{IP})$ & $2$ & $0$ \\
3. & ${\Delta }^{-}_0(0)$ & $1$ & $1$ \\
4. & ${\Delta }^{+}(0)$ & $1$ & $0$ 
\end{tabular}
\end{center}
\begin{center}
Table 1. Possible $\oo (V, K)$ indecomposable types
\end{center}
    
\noindent In table 1 we have use the notation ${\Delta }_0 (\zeta , \mathrm{RP}) = 
{\Delta }_0(\zeta , -\zeta )$, $\zeta = \overline{\zeta }$, $\zeta \ne 0$, where $\zeta $ is a real 
eigenvalue of $Y$ and ${\Delta }_0(\mathrm{i}\beta , \mathrm{IP}) = 
{\Delta }_0(\zeta , - \zeta )$, $\zeta = - \overline{\zeta } = \mathrm{i}\beta $, $\zeta \ne 0$, where 
$\zeta$ is a purely imaginary eigenvalue of $Y$ with $Y^2 - {\beta }^2 =0$. The indecomposable type 
${\Delta }^{\varepsilon }_0(0)$, where ${\varepsilon }^2=1$, 
is represented by the tuple $\big( {\spann }_{\R } \{ f \} , 0; (\varepsilon ) \big) $. Table 1 is 
taken from the classification of indecomposable types in \cite[p.349]{burgoyne-cushman}. In table 2 
we list the possible indecomposable affine types. 
\begin{center}
\begin{tabular}{rccc}
& \multicolumn{1}{c}{type} & \multicolumn{1}{c}{dimension} & \multicolumn{1}{c}{index} \\ \hline
\rule{0pt}{11pt}1. & ${\nabla }^{+}_3(0), \, \mu \ne 0$ & $3$ & $1$ \\
2. & ${\nabla}_2(0,0)$ & $2$ & $1$ 
\end{tabular}
\end{center}
\begin{center}
Table 2. Possible $\oo (V, K)$ indecomposable affine cotypes
\end{center}

\noindent Table 2 is taken from the classification of indecomposable affine cotypes 
in \cite[p.15]{cushman-vanderkallen}. \medskip 

Using the results of \S 4 we get a list of nonzero indecomposable cotypes of ${\mathrm{Gal}}_3$, which gives a 
list of the coadjoint orbits of ${\mathrm{Gal}}_3$. \bigskip

\noindent \hspace{-.5in}\begin{tabular}{rlll}
& \multicolumn{1}{c}{indecomposable cotype} & \multicolumn{1}{c}{dimension} & \multicolumn{1}{c}{index} \\ \hline
\rule{0pt}{11pt}$1$. & ${\nabla}^{+}_3(0), \mu \ne 0 +{\Delta }_0(\mathrm{i}\beta , \mathrm{IP})$ & 
$3+2$ & $1+0$ \\ 
\rule{0pt}{11pt} $2$. & ${\nabla}^{+}_3(0), \mu \ne 0 + {\Delta }^{+}_0(0) + {\Delta }^{+}_0(0)$ & 
$3+1+1$ & $1 + 0+0$ \\
\rule{0pt}{11pt} $3$. & ${\nabla}^y_2(0), y \ne 0  +{\Delta }^{-}_2(0)$ & $2 +3$ & $0 +1$ \\
\rule{0pt}{11pt} $4$. &${\nabla}^y_2(0), y \ne 0  + {\Delta }_0(\zeta , \mathrm{RP}) +{\Delta }^{+}_0(0)$ & 
$2 + 2 +1$ & $0 +1 +0$ \\ 
\rule{0pt}{11pt} $5$. & ${\nabla}^y_2(0), y \ne 0  + {\Delta }_0(\mathrm{i}\beta , \mathrm{IP}) +{\Delta }^{-}_0(0)$ & 
$2 + 2 +1$ & $0+0+1$ \\
\rule{0pt}{11pt} $6$. & ${\nabla }_2(0,0) +{\Delta }_0(\mathrm{i}\beta , \mathrm{IP}) + {\Delta }^{+}_0(0)$ & 
$2 +2 +1$ & $1 +0 +0$ \\
\rule{0pt}{11pt} $7$. & ${\nabla }_2(0,0) + {\Delta }^{+}_0(0) +{\Delta }^{+}_0(0)+ {\Delta }^{+}_0(0)$ & 
$2 +1 +1+1$ & $1+0+0+0$ \\ 
\rule{0pt}{11pt} $8$. &${\nabla}^y_2(0), y \ne 0  + {\Delta }^{-}_0(0)+{\Delta }^{+}_0(0) +{\Delta }^{+}_0(0)$ & 
$2+1+1+1$ & $0 + 1 +0 +0$  \\
\rule{0pt}{11pt} $9$. & ${\nabla}^{-}_0(0) + {\Delta }_0(\mathrm{i}\beta , \mathrm{IP}) + {\Delta }^{+}_0(0) + 
{\Delta }^{+}_0(0)$ & $1 +2 +1 +1$ & $1 + 0 +0 +0$ \\
\rule{0pt}{11pt} $10$. & ${\nabla}^{-}_0(0)+{\Delta }^{+}_0(0) + {\Delta }^{+}_0(0) +{\Delta }^{+}_0(0) +
{\Delta }^{+}_0(0)$ & $1+1+1+1+1$ & $1+0+0+0+0$ \\
\rule{0pt}{11pt} $11$. & ${\nabla}^{+}_0(0)+ {\Delta }^{-}_2(0) + {\Delta }^{+}_0(0)$ & $1+3+1$ & $0+1+0$ \\
\rule{0pt}{11pt} $12$. & ${\nabla}^{+}_0(0) + {\Delta }_0(\zeta , \mathrm{RP}) + 
{\Delta }_0(\mathrm{i}\beta , \mathrm{IP}) $ & $1 +2 +2$ & $0 +1 +0$ \\
\rule{0pt}{11pt} $13$. & ${\nabla}^{+}_0(0) +{\Delta }_0(\zeta , \mathrm{RP}) +{\Delta}^{+}_0(0) + 
{\Delta}^{+}_0(0)$ & $1+2+1+1$ & $0 +1 +0 +0$ \\
\rule{0pt}{11pt} $14$. & ${\nabla}^{+}_0(0) +{\Delta }_0(\mathrm{i}\beta  , \mathrm{IP}) + 
{\Delta}^{-}_0(0) + {\Delta}^{+}_0(0)$ & $1+2+1+1$ & $0 +0 +1 +0$ \\
\rule{0pt}{11pt} $15$. & ${\nabla}^{+}_0(0)+{\Delta}^{-}_0(0) +{\Delta}^{+}_0(0)+{\Delta}^{+}_0(0)+
{\Delta}^{+}_0(0)$ & $1+1+1+1+1$ & $0 +1+0+0+0$ 
\end{tabular} \bigskip 

\hspace{.3in}Table 3. Coadjoint orbits of the general Galile group ${\mathrm{Gal}}_3$.

\end{document}